\documentclass[11pt]{amsart}

\usepackage{amsmath}
\usepackage{amssymb}
\usepackage{amsfonts}
\usepackage{mathrsfs}

\vfuzz2pt \DeclareMathOperator{\RE}{Re}
\DeclareMathOperator{\IM}{Im}

\DeclareMathOperator{\Psl}{PSL} 

 \DeclareMathOperator{\Ad}{Ad}
 \DeclareMathOperator{\sign}{sign}
\DeclareMathOperator{\TRR}{Tr} \DeclareMathOperator{\MOD}{mod}
\DeclareMathOperator{\vol}{vol} 

\DeclareMathOperator{\II}{i}

\newenvironment{teorem}[2][Theorem]{\begin{trivlist}
\item[\hskip \labelsep {\bfseries #1}\hskip \labelsep {\bfseries #2}]}{\end{trivlist}}

\newenvironment{lem}[2][Lemma]{\begin{trivlist}
\item[\hskip \labelsep {\bfseries #1}\hskip \labelsep {\bfseries #2}]}{\end{trivlist}}

\newenvironment{remm}[2][Remark]{\begin{trivlist}
\item[\hskip \labelsep {\bfseries #1}\hskip \labelsep {\bfseries #2}]}{\end{trivlist}}

\numberwithin{equation}{section}
\numberwithin{theorem}{section}

\newcommand{\set}[1]{\left\{#1\right\}}
\newcommand{\abs}[1]{\left\vert#1\right\vert}
\newcommand{\br}[1]{\left(#1\right)}
\newcommand{\SqBr}[1]{\left[#1\right]}

\begin{document}

\title[Logarithmic derivative of the zeta functions]{On the logarithmic derivative of zeta functions for compact even-dimensional locally symmetric spaces}
\author{Muharem Avdispahi\'c and D\v zenan Gu\v si\'c}

\address{University of Sarajevo, Department of Mathematics, Zmaja od Bosne
35, 71000 Sarajevo, Bosnia and Herzegovina}
\email{\textbf{mavdispa@pmf.unsa.ba}}

\address{University of Sarajevo, Department of Mathematics, Zmaja od Bosne
35, 71000 Sarajevo, Bosnia and Herzegovina}
\email{\textbf{dzenang@pmf.unsa.ba}}



\maketitle

\begin{abstract}
 We derive approximate formulas for the logarithmic derivative of the Selberg and Ruelle zeta functions over compact, even-dimensional, locally symmetric spaces of rank one. The obtained formulas are given in terms of the zeta-singularities.
\end{abstract}

\section{Introduction}

Let $Y=\Gamma\backslash G/K=\Gamma\backslash X$ be a compact, $n-$ dimensional ($n$ even), locally
symmetric Riemannian manifold with negative sectional curvature, where $G$ is
a connected semisimple Lie group of real rank one, $K$ is a
maximal compact subgroup of $G$ and $\Gamma$ is a discrete
co-compact torsion free subgroup of $G$. The covering manifold $X$ is known to be a real, a complex or a quaternionic hyperbolic space or the hyperbolic Cayley plane, i.e. $X$ is one of the following spaces:

\[H\mathbb{R}^{k},\,\,H\mathbb{C}^{m},\,\,H\mathbb{H}^{m},\,\,H\mathbb{C}a^{2}.\]
\newline
Here, $n=k$, $2m$, $4m$, $16$, respectively.

We require $G$ to be linear in order to have complexification available.

U. Bunke and M. Olbrich \cite{Bunke} derived the properties of the zeta functions of
Selberg and Ruelle canonically associated with the geodesic flow of $Y$.

In many applications it is often useful to have some approximate representation of the logarithmic derivative of an appropriate zeta function. Following traditional approach \cite{Titchmarsh}, (see also, \cite{Randol2}), we obtain such representations for the zeta functions described in \cite{Bunke}.

\section{Preliminaries}

The notation that will be applied in the sequel follows \cite{Bunke} (see also \cite{AG}, \cite{AG1}).

Let $\mathfrak{g}=\mathfrak{k}\oplus\mathfrak{p}$ be the Cartan
decomposition of the Lie algebra $\mathfrak{g}$ of $G$,
$\mathfrak{a}$ a maximal abelian subspace of $\mathfrak{p}$ and
$M$ the centralizer of $\mathfrak{a}$ in $K$ with Lie algebra
$\mathfrak{m}$. We normalize the $\Ad{\br{G}}-$ invariant inner
product $(.,.)$ on $\mathfrak{g}$ to restrict to the metric on
$\mathfrak{p}$. Let $SX=G/M$ be the unit sphere bundle of $X$.
Hence $SY=\Gamma\backslash G/M$.

Let $\Phi\br{\mathfrak{g},\mathfrak{a}}$ be the root system and
$W=W\br{\mathfrak{g},\mathfrak{a}}\cong\mathbb{Z}_{2}$ its Weyl
group. Fix a system of positive roots
$\Phi^{+}\br{\mathfrak{g},\mathfrak{a}}\subset\Phi\br{\mathfrak{g},\mathfrak{a}}$.
Let

\[\mathfrak{n}=\sum\limits_{\alpha\in\Phi^{+}\br{\mathfrak{g},\mathfrak{a}}}\mathfrak{n}_{\alpha}\]
\newline
be the sum of the root spaces corresponding to elements of
$\Phi^{+}\br{\mathfrak{g},\mathfrak{a}}$. The decomposition
$\mathfrak{g}=\mathfrak{k}\oplus\mathfrak{a}\oplus\mathfrak{n}$
corresponds to the Iwasawa decomposition of the group $G=KAN$.
Define $\rho\in\mathfrak{a}_{\mathbb{C}}^{*}$ by

\[\rho=\frac{1}{2}\sum\limits_{\alpha\in\Phi^{+}\br{\mathfrak{g},\mathfrak{a}}}\dim\br{\mathfrak{n}_{\alpha}}\alpha.\]
\newline

We normalize the metric on $Y$ to be of sectional curvature $-1$ if $X=H\mathbb{R}^{k}$. In all other cases, the normalized metric on $Y$ is such that the sectional curvature varies between $-1$ and $-4$. Hence, $\rho=\frac{1}{2}\br{k-1}$, $m$, $2m+1$, $11$ if $n=k$, $2m$, $4m$, $16$, respectively.

The positive Weyl chamber $\mathfrak{a}^{+}$ is the half line in
$\mathfrak{a}$ on which the positive roots take positive values.
Let $A^{+}=\exp\br{\mathfrak{a}^{+}}\subset A$.

The symmetric space $X$ has a compact dual space $X_{d}=G_{d}/K$,
where $G_{d}$ is the analytic subgroup of $GL\br{n,\mathbb{C}}$
corresponding to
$\mathfrak{g}_{d}=\mathfrak{k}\oplus\mathfrak{p}_{d}$,
$\mathfrak{p}_{d}=$i$\mathfrak{p}$. We normalize the metric on
$X_{d}$ in such a way that the multiplication by i induces an
isometry between $\mathfrak{p}$ and $\mathfrak{p}_{d}$.

Let $i^{*}:R\br{K}\rightarrow R\br{M}$ be the restriction map
induced by the embedding $i:M\hookrightarrow K$, where $R\br{K}$
and $R\br{M}$ are the representation rings over $\mathbb{Z}$ of
$K$ and $M$, respectively.

Since $n$ is even, every $\sigma\in\hat{M}$ is invariant under the
action of the Weyl group $W$ (see, \cite[p.~27]{Bunke}). Let
$\sigma\in\hat{M}$. We choose $\gamma\in R\br{K}$ such that
$i^{*}\br{\gamma}=\sigma$ and represent it by $\Sigma
a_{i}\gamma_{i}$, $a_{i}\in\mathbb{Z},\gamma_{i}\in\hat{K}$. Set

\[V_{\gamma}^{\pm}=\sum\limits_{\sign{\br{a_{i}}}=\pm1}\sum\limits_{m=1}^{\abs{a_{i}}}V_{\gamma_{i}},\]
\newline
where $V_{\gamma_{i}}$ is the representation space of
$\gamma_{i}$. Define
$V\br{\gamma}^{\pm}=G\times_{K}V_{\gamma}^{\pm}$ and
$V_{d}\br{\gamma}^{\pm}=G_{d}\times_{K}V_{\gamma}^{\pm}$. To
$\gamma$ we associate $\mathbb{Z}_{2}-$ graded homogeneous vector
bundles $V\br{\gamma}=V\br{\gamma}^{+}\oplus V\br{\gamma}^{-}$ and
$V_{d}\br{\gamma}=V_{d}\br{\gamma}^{+}\oplus V_{d}\br{\gamma}^{-}$
on $X$ and $X_{d}$, respectively. Let
\[V_{Y,\chi}\br{\gamma}=\Gamma\backslash\br{V_{\chi}\otimes
V\br{\gamma}}\]
\newline
be a $\mathbb{Z}_{2}-$ graded vector bundle on $Y$, where
$\br{\chi,V_{\chi}}$ is a finite-dimensional unitary
representation of $\Gamma$.

Reasoning as in \cite[beginning of Subsection 1.1.2]{Bunke},
we choose a Cartan subalgebra $\mathfrak{t}$ of $\mathfrak{m}$ and
a system of positive roots
$\Phi^{+}\br{\mathfrak{m}_{\mathbb{C}},\mathfrak{t}}$. Then,
$\rho_{\mathfrak{m}}\in$i$\mathfrak{t}^{*}$, where
\[\rho_{\mathfrak{m}}=\frac{1}{2}\sum\limits_{\alpha\in\Phi^{+}\br{\mathfrak{m}_{\mathbb{C}},\mathfrak{t}}}\alpha.\]
\newline
Let $\mu_{\sigma}\in$i$\mathfrak{t}^{*}$ be the highest weight of
$\sigma$. Set

\[c\br{\sigma}=\abs{\rho}^{2}+\abs{\rho_{\mathfrak{m}}}^{2}-\abs{\mu_{\sigma}+\rho_{\mathfrak{m}}}^{2},\]
\newline
where the norms are induced by the complex bilinear extension to
$\mathfrak{g}_{\mathbb{C}}$ of the inner product $(.,.)$. Finally,
we introduce the operators  (see, \cite[p.~28]{Bunke})

\[A_{d}\br{\gamma,\sigma}^{2}=\Omega +c\br{\sigma}: C^{\infty}\br{X_{d},V_{d}\br{\gamma}}\rightarrow C^{\infty}\br{X_{d},V_{d}\br{\gamma}},\]
\[A_{Y,\chi}\br{\gamma,\sigma}^{2}=-\Omega -c\br{\sigma}: C^{\infty}\br{Y,V_{Y,\chi}\br{\gamma}}\rightarrow C^{\infty}\br{Y,V_{Y,\chi}\br{\gamma}},\]
\newline
$\Omega$ being the Casimir element of the complex
universal enveloping algebra $\mathcal{U}\br{\mathfrak{g}}$ of
$\mathfrak{g}$.

Let $m_{\chi}\br{s,\gamma,\sigma}=\TRR{E_{A_{Y,\chi}\br{\gamma,\sigma}}}\br{\set{s}}$,
$m_{d}\br{s,\gamma,\sigma}=\TRR{E_{A_{d}\br{\gamma,\sigma}}}\br{\set{s}}$,
where $E_{A}\br{.}$ denotes the family of spectral projections of
a normal operator $A$.

Now, we choose a maximal abelian subalgebra $\mathfrak{t}$ of
$\mathfrak{m}$. Then,
$\mathfrak{h}=\mathfrak{t}_{\mathbb{C}}\oplus\mathfrak{a}_{\mathbb{C}}$
is a Cartan subalgebra of $\mathfrak{g}_{\mathbb{C}}$. Let
$\Phi^{+}\br{\mathfrak{g}_{\mathbb{C}},\mathfrak{h}}$ be a
positive root system having the property that, for
$\alpha\in\Phi\br{\mathfrak{g}_{\mathbb{C}},\mathfrak{h}}$,
$\alpha_{|\mathfrak{a}}\in\Phi^{+}\br{\mathfrak{g},\mathfrak{a}}$
implies
$\alpha\in\Phi^{+}\br{\mathfrak{g}_{\mathbb{C}},\mathfrak{h}}$.
Let
\[\delta=\frac{1}{2}\sum\limits_{\alpha\in\Phi^{+}\br{\mathfrak{g}_{\mathbb{C}},\mathfrak{h}}}\alpha.\]
\newline
We set $\rho_{\mathfrak{m}}=\delta-\rho$. Define the root vector
$H_{\alpha}\in\mathfrak{a}$ for
$\alpha\in\Phi^{+}\br{\mathfrak{g},\mathfrak{a}}$ by

\[\lambda\br{H_{\alpha}}=\frac{\br{\lambda,\alpha}}{\br{\alpha,\alpha}},\]
where $\lambda\in\mathfrak{a}^{*}$.

For $\alpha\in\Phi^{+}\br{\mathfrak{g},\mathfrak{a}}$, we define
$\varepsilon_{\alpha}\br{\sigma}\in\set{0,\frac{1}{2}}$ by

\[e^{2\pi \textrm{i}\varepsilon_{\alpha}\br{\sigma}}=\sigma\br{e^{2\pi
\textrm{i}H_{\alpha}}}\in\set{\pm 1}.\]
\newline
According to \cite[p.~47]{Bunke}, the root system
$\Phi^{+}\br{\mathfrak{g},\mathfrak{a}}$ is of the form
$\Phi^{+}\br{\mathfrak{g},\mathfrak{a}}=\set{\alpha}$ or
$\Phi^{+}\br{\mathfrak{g},\mathfrak{a}}=\set{\frac{\alpha}{2},\alpha}$
for the long root $\alpha$. Let $\alpha$ be the long root in
$\Phi^{+}\br{\mathfrak{g},\mathfrak{a}}$. We set $T=\abs{\alpha}$.
For $\sigma\in\hat{M}$, $\epsilon_{\sigma}\in\set{0,\frac{1}{2}}$
is given by

\[\epsilon_{\sigma}\equiv\frac{\abs{\rho}}{T}+\varepsilon_{\alpha}\br{\sigma}\,\MOD{\mathbb{Z}}.\]
\newline
We define the lattice
$L\br{\sigma}\subset\mathbb{R}\cong\mathfrak{a}^{*}$ by
$L\br{\sigma}=T\br{\epsilon_{\sigma}+\mathbb{Z}}$. Finally, for
$\lambda\in\mathfrak{a}_{\mathbb{C}}^{*}\cong\mathbb{C}$ we set
\[P_{\sigma}\br{\lambda}=\prod\limits_{\beta\in\Phi^{+}\br{\mathfrak{g}_{\mathbb{C}},\mathfrak{h}}}
\frac{\br{\lambda
+\mu_{\sigma}+\rho_{\mathfrak{m}},\beta}}{\br{\delta,\beta}}.\]
\newline

Since $n$ is even, there exists a $\sigma-$ admissible $\gamma\in
R\br{K}$ for every $\sigma\in\hat{M}$ (see, \cite[p.~49, Lemma
1.18]{Bunke}). Here, $\gamma\in R\br{K}$ is called $\sigma-$
admissible if $i^{*}\br{\gamma}=\sigma$ and
$m_{d}\br{s,\gamma,\sigma}=P_{\sigma}\br{s}$ for all $0\leq s\in
L\br{\sigma}$.

\section{Zeta functions and the geodesic flow}

Since $\Gamma\subset G$ is co-compact and torsion free, there are
only two types of conjugacy classes - the class of the identity
$1\in\Gamma$ and classes of hyperbolic elements.

Let $g\in G$ be hyperbolic. Then there is an Iwasawa decomposition
$G=NAK$ such that $g=am\in A^{+}M$. Following \cite[p.~59]{Bunke},
we define

\[l\br{g}=\abs{\log\br{a}}.\]
\newline
Let $\Gamma_{\textrm{h}}$, resp. $\text{P}\Gamma_{\textrm{h}}$
denote the set of the $\Gamma -$ conjugacy classes of hyperbolic
resp. primitive hyperbolic elements in $\Gamma$.

Let $\varphi$ be the geodesic flow on $SY$ determined by the
metric of $Y$. In the representation $SY=\Gamma\backslash G/M$, $\varphi$
is given by

\[\varphi : \mathbb{R}\times SY\ni\br{t,\Gamma gM}\rightarrow\Gamma g\exp\br{-tH}M\in SY,\]
\newline
where $H$ is the unit vector in $\mathfrak{a}^{+}$. If
$V_{\chi}\br{\sigma}=\Gamma\backslash\br{G\times_{M}V_{\sigma}\otimes
V_{\chi}}$ is the vector bundle corresponding to
finite-dimensional unitary representations
$\br{\sigma,V_{\sigma}}$ of $M$ and $\br{\chi,V_{\chi}}$ of
$\Gamma$, then we define a lift $\varphi_{\chi,\sigma}$ of
$\varphi$ to $V_{\chi}\br{\sigma}$ by (see, \cite[p.~95]{Bunke})

\[\varphi_{\chi,\sigma} : \mathbb{R}\times V_{\chi}\br{\sigma}\ni\br{t,\SqBr{g,v\otimes w}}\rightarrow
\SqBr{g\exp\br{-tH},v\otimes w}\in V_{\chi}\br{\sigma}.\]
\newline
For $\RE{\br{s}}>2\rho$, the Ruelle zeta function for the flow
$\varphi_{\chi,\sigma}$ is defined by the infinite product

\[Z_{R,\chi}\br{s,\sigma}=\prod\limits_{\gamma_{0}\in\text{P}\Gamma_{\textrm{h}}}
\det\br{1-\br{\sigma\br{m}\otimes\chi\br{\gamma_{0}}}e^{-sl\br{\gamma_{0}}}}^{\br{-1}^{n-1}}.\]
\newline
The Selberg zeta function for the flow $\varphi_{\chi,\sigma}$ is
given by

\[Z_{S,\chi}\br{s,\sigma}=\]
\[\prod\limits_{\gamma_{0}\in\text{P}\Gamma_{\textrm{h}}}\prod\limits_{k=0}^{+\infty}
\det\br{1-\br{\sigma\br{m}\otimes\chi\br{\gamma_{0}}\otimes
S^{k}\br{\Ad{\br{ma}_{\bar{\mathfrak{n}}}}}}e^{-\br{s+\rho}l\br{\gamma_{0}}}},\]
\newline
for $\RE{\br{s}}>\rho$, where $S^{k}$ denotes the $k-$th symmetric
power of an endomorphism, $\bar{\mathfrak{n}}=\theta\mathfrak{n}$
is the sum of negative root spaces of $\mathfrak{a}$ as usual, and
$\theta$ is the Cartan involution of $\mathfrak{g}$.

Let $\mathfrak{n}_{\mathbb{C}}$ be the complexification of
$\mathfrak{n}$. For
$\lambda\in\mathbb{C}\cong\mathfrak{a}_{\mathbb{C}}^{*}$ let
$\mathbb{C}_{\lambda}$ denote the one-dimensional representation
of $A$ given by $A\ni a\rightarrow a^{\lambda}$. Let $p\geq 0$.
There exist sets
\[I_{p}=\set{\br{\tau,\lambda}\mid\tau\in\hat{M},\lambda\in\mathbb{R}}\]
\newline
such that $\Lambda^{p}\mathfrak{n}_{\mathbb{C}}$ as a
representation of $MA$ decomposes with respect to $MA$ as
\[\Lambda^{p}\mathfrak{n}_{\mathbb{C}}=\sum\limits_{\br{\tau,\lambda}\in I_{p}}V_{\tau}\otimes\mathbb{C}_{\lambda},\]
\newline
where $V_{\tau}$ is the space of the representation $\tau$. Bunke
and Olbrich proved that the Ruelle zeta function $Z_{R,\chi}\br{s,\sigma}$ has the following
representation (see, \cite[p.~99, Prop. 3.4]{Bunke})

\begin{equation}\label{3.1}
Z_{R,\chi}\br{s,\sigma}=\prod\limits_{p=0}^{n-1}\prod\limits_{\br{\tau,\lambda}\in I_{p}}Z_{S,\chi}\br{s+\rho-\lambda,\tau\otimes\sigma}^{\br{-1}^{p}}.
\end{equation}

Let $d_{Y}=-\br{-1}^{\frac{n}{2}}$. The following theorem holds
true (see, \cite[p.~113, Th. 3.15]{Bunke}).
\newline
\begin{teorem}{A.}\label{t3.15}
\textit{The Selberg zeta function $Z_{S,\chi}\br{s,\sigma}$ has a
meromorphic continuation to all of $\mathbb{C}$. If $\gamma$ is
$\sigma-$admissible, then the singularities (zeros and poles) of
$Z_{S,\chi}\br{s,\sigma}$ are the following ones:\\
\begin{enumerate}
    \item at $\pm$ $\II$$s$ of order $m_{\chi}\br{s,\gamma,\sigma}$ if
    $s\neq 0$ is an eigenvalue of
    $A_{Y,\chi}\br{\gamma,\sigma}$,\\
    \item at $s=0$ of order $2m_{\chi}\br{0,\gamma,\sigma}$ if
    $0$ is an eigenvalue of $A_{Y,\chi}\br{\gamma,\sigma}$,\\
    \item at $-s$, $s\in T\br{\mathbb{N}-\epsilon_{\sigma}}$ of
    order
    $2\frac{d_{Y}\dim\br{\chi}\vol\br{Y}}{\vol\br{X_{d}}}m_{d}\br{s,\gamma,\sigma}$. Then $s>0$ is an eigenvalue of $A_{d}\br{\gamma,\sigma}$.\\
\end{enumerate}
If two such points coincide, then the orders add up.}
\end{teorem}\

Note that the shifts $\rho-\lambda$ that appear in (\ref{3.1}) are always contained in the interval $\SqBr{-\rho,\rho}$, (see, \cite{AG1}).

In \cite{AG}, we proved that there exist entire functions $Z_{S}^{1}\br{s}$, $Z_{S}^{2}\br{s}$ of order at most $n$ such that
\begin{equation}\label{3.2}
Z_{S,\chi}\br{s,\sigma}=\frac{Z_{S}^{1}\br{s}}{Z_{S}^{2}\br{s}}.
\end{equation}\
\newline
Here, $\gamma$ is $\sigma\,-$ admissible, the zeros of $Z_{S}^{1}\br{s}$ correspond to the
zeros of $Z_{S,\chi}\br{s,\sigma}$ and the zeros of
$Z_{S}^{2}\br{s}$ correspond to the poles of
$Z_{S,\chi}\br{s,\sigma}$. The orders of the zeros of
$Z_{S}^{1}\br{s}$ resp. $Z_{S}^{2}\br{s}$ equal the orders of the
corresponding zeros resp. poles of $Z_{S,\chi}\br{s,\sigma}$. Furthermore, (see, \cite{AG1}),

\begin{equation}\label{3.3}
\abs{Z^{i}_{S}\br{\sigma_{1}+\II t}}=e^{O\br{\abs{t}^{n-1}}}
\end{equation}\
\newline
uniformly in any bounded strip $b_{1}\leq\sigma_{1}\leq b_{2}$ for $i=1,2$.

Denote by $N\br{t}$ the number of singularities of $Z_{S,\chi}\br{s,\sigma}$ on the interval $\II$$x$, $0<x<t$. In \cite{AG1}, we proved that

\begin{equation}\label{3.4}
N\br{t}=\frac{\dim\br{\chi}\vol{\br{Y}}}{nT\vol{\br{X_{d}}}}t^{n}+O\br{t^{n-1}}.
\end{equation}\
\newline
Moreover, we proved that there exists a constant $C$ such that

\begin{equation}\label{3.5}
N_{R}\br{t}=Ct^{n}+O\br{t^{n-1}},
\end{equation}\
\newline
where $N_{R}\br{t}$ denotes the number of singularities of $Z_{R,\chi}\br{s,\sigma}$ in the rectangle $a\leq\RE{\br{s}}\leq b$, $0<\IM{\br{s}}<t$. Here, $-\rho\leq a\leq b\leq\rho$.

Recall that $\gamma$ is assumed to be $\sigma\,-$ admissible in (\ref{3.4}) and (\ref{3.5}).

The following well known lemma will be used in the sequel (see \cite[p.~56]{Titchmarsh})
\newline
\begin{lem}{B.}\label{lll}
\textit{If $f\br{s}$ is regular, and}

\[\abs{\frac{f\br{s}}{f\br{s_{0}}}}<e^{M}\,\,\br{M>1}\]
\newline
\textit{in the circle $\abs{s-s_{0}}\leq r$, then}

\[\abs{\frac{f^{'}\br{s}}{f\br{s}}-\sum\limits_{\rho}\frac{1}{s-\rho}}<\frac{AM}{r}\,\,\br{\abs{s-s_{0}}\leq\frac{1}{4}r},\]
\newline
\textit{where $\rho$ runs through the zeros of $f\br{s}$ such that $\abs{\rho-s_{0}}\leq\frac{1}{2}r$}.
\end{lem}

\section{Main result}

The main result of the paper is the following theorem.
\newline
\begin{teorem}{4.1.}\label{MT}
\textit{Let $\gamma$ be $\sigma\,-$ admissible. Suppose $t\gg 0$ is selected so that $\rho+\II t$ is not a singularity of $Z_{R,\chi}\br{s,\sigma}$. Then,}
\begin{enumerate}
\item[(\textit{a})]
\[\frac{Z^{'}_{R,\chi}\br{s,\sigma}}{Z_{R,\chi}\br{s,\sigma}}=O\br{t^{n-1}}+
\sum_{\abs{t-t_{R}}\leq 1}\frac{1}{s-s_{R}},\]
\newline
\textit{where $s=\sigma_{1}+\II t$, $\rho\leq\sigma_{1}\leq\frac{1}{2}t+\rho$ and $s_{R}=\rho+\II t_{R}$ is a singularity of $Z_{R,\chi}\br{s,\sigma}$ along the line $\RE{\br{s}}=\rho$.}
\item[(\textit{b})]
\[\frac{Z^{'}_{R,\chi}\br{s,\sigma}}{Z_{R,\chi}\br{s,\sigma}}=O\br{t^{n-1}},\]
\newline
\textit{where $s=\sigma_{1}+\II t$, $\rho+u\leq\sigma_{1}\leq\frac{1}{2}t+\rho$ and $u>0$.}
    \end{enumerate}
\end{teorem}\
\begin{proof}
(\textit{a}) Let $\br{\tau,2\rho}\in I_{p}$ for some $p\in\set{0,1,...,n-1}$. Then,

\[Z_{S,\chi}\br{s-\rho,\tau\otimes\sigma}^{\br{-1}^{p}}\]
\newline
is the corresponding factor in the representation (\ref{3.1}). By (\ref{3.2}) and (\ref{3.3}),

\begin{equation}\label{4.1}
\abs{Z_{S,\chi}\br{s-\rho,\tau\otimes\sigma}}=e^{O\br{t_{1}^{n-1}}}
\end{equation}
\newline
uniformly in any bounded half-strip $b_{1}\leq\RE{\br{s}}\leq b_{2}$, $s=\sigma_{1}+\II t_{1}$, $t_{1}>0$.

Let $8\rho\leq r<t$. We choose $c$, $2\rho<c<\frac{1}{4}r+\rho$ and put $s_{0}=c+\II t$. It follows immediately that the circles $\abs{s-s_{0}}\leq r$, $\abs{s-s_{0}}\leq\frac{1}{2}r$ and $\abs{s-s_{0}}\leq\frac{1}{4}r$ cross the line $\RE{\br{s}}=\rho$.

Denote the set of poles of $Z_{S,\chi}\br{s-\rho,\tau\otimes\sigma}$ lying in the circle $\abs{s-s_{0}}\leq r$ by $P$. Then, the function

\[\mathcal{H}\br{s}=Z_{S,\chi}\br{s-\rho,\tau\otimes\sigma}\cdot\prod\limits_{\rho_{1}\in P}\br{s-\rho_{1}}\]
\newline
is regular in $\abs{s-s_{0}}\leq r$. By (\ref{4.1}),

\[\abs{Z_{S,\chi}\br{s-\rho,\tau\otimes\sigma}}=e^{O\br{t_{1}^{n-1}}}\]
\newline
uniformly in the half-strip $c-r\leq\RE{\br{s}}\leq c+r$, $s=\sigma_{1}+\II t_{1}$, $t_{1}>0$. Hence,

\begin{equation}\label{4.2}
\abs{Z_{S,\chi}\br{s-\rho,\tau\otimes\sigma}}=e^{O\br{t_{1}^{n-1}}}
\end{equation}
\newline
for $s=\sigma_{1}+\II t_{1}$, $\abs{s-s_{0}}\leq r$. Specially,

\begin{equation}\label{4.3}
\abs{Z_{S,\chi}\br{s_{0}-\rho,\tau\otimes\sigma}}=e^{O\br{t^{n-1}}}.
\end{equation}
\newline
Having in mind that $t_{1}\leq t+r<2t$ for $s=\sigma_{1}+\II t_{1}$, $\abs{s-s_{0}}\leq r$, the relations (\ref{4.2}) and (\ref{4.3}) imply

\begin{equation}\label{4.4}
\abs{\frac{Z_{S,\chi}\br{s-\rho,\tau\otimes\sigma}}{Z_{S,\chi}\br{s_{0}-\rho,\tau\otimes\sigma}}}=e^{O\br{t^{n-1}}}
\end{equation}
\newline
for $s=\sigma_{1}+\II t_{1}$, $\abs{s-s_{0}}\leq r$.

Since $P$ is a finite set and $\abs{s-\rho_{1}}\leq 2r$, $\abs{s_{0}-\rho_{1}}>\rho$ for all $\rho_{1}\in P$ and $s=\sigma_{1}+\II t_{1}$, $\abs{s-s_{0}}\leq r$, it follows from (\ref{4.4}) that
\newline

\[\abs{\frac{\mathcal{H}\br{s}}{\mathcal{H}\br{s_{0}}}}=
\abs{\frac{Z_{S,\chi}\br{s-\rho,\tau\otimes\sigma}}{Z_{S,\chi}\br{s_{0}-\rho,\tau\otimes\sigma}}}\cdot\prod\limits_{\rho_{1}\in P}\frac{\abs{s-\rho_{1}}}{\abs{s_{0}-\rho_{1}}}=\]

\[e^{O\br{t^{n-1}}}\cdot O\br{1}=e^{O\br{t^{n-1}}}\]
\newline
for $s=\sigma_{1}+\II t_{1}$, $\abs{s-s_{0}}\leq r$. Hence, there exists a constant $C$ such that

\[\abs{\frac{\mathcal{H}\br{s}}{\mathcal{H}\br{s_{0}}}}<e^{Ct^{n-1}}\]
\newline
for $s=\sigma_{1}+\II t_{1}$, $\abs{s-s_{0}}\leq r$. Putting $M=Ct^{n-1}$ and applying Lemma B, we obtain
\[\frac{\mathcal{H}^{'}\br{s}}{\mathcal{H}\br{s}}=O\br{t^{n-1}}+\sum\limits_{\rho_{2}\in Q}\frac{1}{s-\rho_{2}}\]
\newline
for $s=\sigma_{1}+\II t_{1}$, $\abs{s-s_{0}}\leq\frac{1}{4}r$, where $Q$ denotes the set of zeros of $\mathcal{H}\br{s}$ lying in $\abs{s-s_{0}}\leq\frac{1}{2}r$. It follows from the definition of $\mathcal{H}\br{s}$  that

\begin{equation}\label{4.5}
\frac{Z^{'}_{S,\chi}\br{s-\rho,\tau\otimes\sigma}}{Z_{S,\chi}\br{s-\rho,\tau\otimes\sigma}}=
O\br{t^{n-1}}+\sum\limits_{\rho_{2}\in Q}\frac{1}{s-\rho_{2}}-\sum\limits_{\rho_{1}\in P}\frac{1}{s-\rho_{1}}
\end{equation}
\newline
for $s=\sigma_{1}+\II t_{1}$, $\abs{s-s_{0}}\leq\frac{1}{4}r$. In particular, (\ref{4.5}) remains valid for $s=\sigma_{1}+\II t$, $\rho\leq\sigma_{1}<c+\frac{1}{4}r$. However, $c<\frac{1}{4}r+\rho$. Hence, we get

\begin{equation}\label{4.6}
\frac{Z^{'}_{S,\chi}\br{s-\rho,\tau\otimes\sigma}}{Z_{S,\chi}\br{s-\rho,\tau\otimes\sigma}}=
O\br{t^{n-1}}+\sum\limits_{\rho_{2}\in Q}\frac{1}{s-\rho_{2}}-\sum\limits_{\rho_{1}\in P}\frac{1}{s-\rho_{1}}
\end{equation}
\newline
for $s=\sigma_{1}+\II t$, $\rho\leq\sigma_{1}<\frac{1}{2}r+\rho$.

One can see from the definition of $\mathcal{H}\br{s}$ that $Q$ is the set of zeros of $Z_{S,\chi}\br{s-\rho,\tau\otimes\sigma}$ lying in $\abs{s-s_{0}}\leq\frac{1}{2}r$. Put $\rho_{2}=\rho+\II\gamma_{1}$. We have

\[\abs{\rho_{2}-s_{0}}\leq\frac{1}{2}r\]
if and only if

\[t-\sqrt{\frac{1}{4}r^{2}-\br{c-\rho}^{2}}\leq\gamma_{1}\leq t+\sqrt{\frac{1}{4}r^{2}-\br{c-\rho}^{2}}.\]
Note that
\[2\rho<c<\frac{1}{4}r+\rho\]
if and only if
\[\frac{\sqrt{3}}{4}r<\sqrt{\frac{1}{4}r^{2}-\br{c-\rho}^{2}}<\sqrt{\frac{1}{4}r^{2}-\rho^{2}}.\]
\newline
Taking into account our normalization of the metric on $Y$, we obtain

\[\sqrt{\frac{1}{4}r^{2}-\br{c-\rho}^{2}}>\frac{\sqrt{3}}{4}r\geq 2\sqrt{3}\rho>1.\]
\newline
Hence, the first sum on the right hand side of (\ref{4.6}) can be written as

\begin{equation}\label{4.7}
\sum\limits_{\rho_{2}\in Q}\frac{1}{s-\rho_{2}}=\sum\limits_{\abs{t-\gamma_{1}}\leq 1}\frac{1}{s-\rho_{2}}+
\end{equation}
\[\sum\limits_{t+1<\gamma_{1}\leq t+\sqrt{\frac{1}{4}r^{2}-\br{c-\rho}^{2}}}\frac{1}{s-\rho_{2}}+\sum\limits_{t-\sqrt{\frac{1}{4}r^{2}-\br{c-\rho}^{2}}\leq\gamma_{1}< t-1}\frac{1}{s-\rho_{2}}\]
\newline
for $s=\sigma_{1}+\II t$, $\rho\leq\sigma_{1}<\frac{1}{2}r+\rho$.

Similarly,
\[\abs{\rho_{1}-s_{0}}\leq r\]
if and only if
\[t-\sqrt{r^{2}-\br{c-\rho}^{2}}\leq\gamma_{2}\leq t+\sqrt{r^{2}-\br{c-\rho}^{2}},\]
\newline
where $\rho_{1}=\rho+\II\gamma_{2}$. Therefore,

\begin{equation}\label{4.8}
\sum\limits_{\rho_{1}\in P}\frac{1}{s-\rho_{1}}=\sum\limits_{\abs{t-\gamma_{2}}\leq 1}\frac{1}{s-\rho_{1}}+
\end{equation}
\[\sum\limits_{t+1<\gamma_{2}\leq t+\sqrt{r^{2}-\br{c-\rho}^{2}}}\frac{1}{s-\rho_{1}}+\sum\limits_{t-\sqrt{r^{2}-\br{c-\rho}^{2}}\leq\gamma_{2}< t-1}\frac{1}{s-\rho_{1}}\]
\newline
for $s=\sigma_{1}+\II t$, $\rho\leq\sigma_{1}<\frac{1}{2}r+\rho$. Combining (\ref{4.6}), (\ref{4.7}) and (\ref{4.8}), we obtain

\begin{equation}\label{4.9}
\frac{Z^{'}_{S,\chi}\br{s-\rho,\tau\otimes\sigma}}{Z_{S,\chi}\br{s-\rho,\tau\otimes\sigma}}=O\br{t^{n-1}}+\sum\limits_{\abs{t-\gamma_{1}}\leq 1}\frac{1}{s-\rho_{2}}-\sum\limits_{\abs{t-\gamma_{2}}\leq 1}\frac{1}{s-\rho_{1}}
\end{equation}
\[\sum\limits_{t+1<\gamma_{1}\leq t+\sqrt{\frac{1}{4}r^{2}-\br{c-\rho}^{2}}}\frac{1}{s-\rho_{2}}-\sum\limits_{t+1<\gamma_{2}\leq t+\sqrt{r^{2}-\br{c-\rho}^{2}}}\frac{1}{s-\rho_{1}}+\]
\[\sum\limits_{t-\sqrt{\frac{1}{4}r^{2}-\br{c-\rho}^{2}}\leq\gamma_{1}< t-1}\frac{1}{s-\rho_{2}}-\sum\limits_{t-\sqrt{r^{2}-\br{c-\rho}^{2}}\leq\gamma_{2}< t-1}\frac{1}{s-\rho_{1}}\]
\newline
for $s=\sigma_{1}+\II t$, $\rho\leq\sigma_{1}<\frac{1}{2}r+\rho$.

Corresponding to the pair $\br{\tau,2\rho}\in I_{p}$, singularities of $Z_{S,\chi}\br{s-\rho,\tau\otimes\sigma}$ along the line $\RE{\br{s}}=\rho$ resp. the number of singularities on the interval $\rho+\II x$, $0<x<t$  will be denoted by $\rho_{S,p,\tau}=\rho+\II\gamma_{S,p,\tau}$ resp. $N_{S,p,\tau}\br{t}$.

We have

\[\abs{\sum\limits_{t+1<\gamma_{1}\leq t+\sqrt{\frac{1}{4}r^{2}-\br{c-\rho}^{2}}}\frac{1}{s-\rho_{2}}-\sum\limits_{t+1<\gamma_{2}\leq t+\sqrt{r^{2}-\br{c-\rho}^{2}}}\frac{1}{s-\rho_{1}}}\leq\]
\[\sum\limits_{t+1<\gamma_{1}\leq t+\sqrt{\frac{1}{4}r^{2}-\br{c-\rho}^{2}}}\frac{1}{\abs{s-\rho_{2}}}+\sum\limits_{t+1<\gamma_{2}\leq t+\sqrt{r^{2}-\br{c-\rho}^{2}}}\frac{1}{\abs{s-\rho_{1}}}<\]
\[\sum\limits_{t+1<\gamma_{1}\leq t+\sqrt{\frac{1}{4}r^{2}-\br{c-\rho}^{2}}}1+\sum\limits_{t+1<\gamma_{2}\leq t+\sqrt{r^{2}-\br{c-\rho}^{2}}}1<\]
\[\sum\limits_{t+1<\gamma_{1}\leq t+\sqrt{r^{2}-\br{c-\rho}^{2}}}1+\sum\limits_{t+1<\gamma_{2}\leq t+\sqrt{r^{2}-\br{c-\rho}^{2}}}1=\]
\[N_{S,p,\tau}\br{t+\sqrt{r^{2}-\br{c-\rho}^{2}}}-N_{S,p,\tau}\br{t+1}\]
\newline
for $s=\sigma_{1}+\II t$, $\rho\leq\sigma_{1}<\frac{1}{2}r+\rho$. Hence, it follows from (\ref{3.4}) that

\begin{equation}\label{4.10}
\sum\limits_{t+1<\gamma_{1}\leq t+\sqrt{\frac{1}{4}r^{2}-\br{c-\rho}^{2}}}\frac{1}{s-\rho_{2}}-\sum\limits_{t+1<\gamma_{2}\leq t+\sqrt{r^{2}-\br{c-\rho}^{2}}}\frac{1}{s-\rho_{1}}=
\end{equation}
\[O\br{t^{n-1}}\]
\newline
for $s=\sigma_{1}+\II t$, $\rho\leq\sigma_{1}<\frac{1}{2}r+\rho$. Similarly,

\begin{equation}\label{4.11}
\sum\limits_{t-\sqrt{\frac{1}{4}r^{2}-\br{c-\rho}^{2}}\leq\gamma_{1}< t-1}\frac{1}{s-\rho_{2}}-\sum\limits_{t-\sqrt{r^{2}-\br{c-\rho}^{2}}\leq\gamma_{2}< t-1}\frac{1}{s-\rho_{1}}=
\end{equation}
\[O\br{t^{n-1}}\]
\newline
for $s=\sigma_{1}+\II t$, $\rho\leq\sigma_{1}<\frac{1}{2}r+\rho$. Combining (\ref{4.9}), (\ref{4.10}) and (\ref{4.11}), we conclude
\[\frac{Z^{'}_{S,\chi}\br{s-\rho,\tau\otimes\sigma}}{Z_{S,\chi}\br{s-\rho,\tau\otimes\sigma}}=O\br{t^{n-1}}+\sum\limits_{\abs{t-\gamma_{S,p,\tau}}\leq 1}\frac{1}{s-\rho_{S,p,\tau}}\]
\newline
for $s=\sigma_{1}+\II t$, $\rho\leq\sigma_{1}<\frac{1}{2}r+\rho$. However, $r<t$. Hence,

\begin{equation}\label{4.12}
\frac{Z^{'}_{S,\chi}\br{s-\rho,\tau\otimes\sigma}}{Z_{S,\chi}\br{s-\rho,\tau\otimes\sigma}}=O\br{t^{n-1}}+\sum\limits_{\abs{t-\gamma_{S,p,\tau}}\leq 1}\frac{1}{s-\rho_{S,p,\tau}}
\end{equation}
\newline
for $s=\sigma_{1}+\II t$, $\rho\leq\sigma_{1}<\frac{1}{2}t+\rho$.

Let $u>0$. One has

\[\abs{\sum\limits_{\abs{t-\gamma_{S,p,\tau}}\leq 1}\frac{1}{s-\rho_{S,p,\tau}}}\leq\sum\limits_{\abs{t-\gamma_{S,p,\tau}}\leq 1}\frac{1}{\abs{s-\rho_{S,p,\tau}}}<\frac{1}{u}\sum\limits_{\abs{t-\gamma_{S,p,\tau}}\leq 1}1=\]
\[\frac{1}{u}\br{N_{S,p,\tau}\br{t+1}-N_{S,p,\tau}\br{t-1}}\]
\newline
for $s=\sigma_{1}+\II t$, $\rho+u\leq\sigma_{1}<\frac{1}{2}t+\rho$. Therefore, it follows from (\ref{3.4}) and (\ref{4.12}) that
\begin{equation}\label{4.13}
\frac{Z^{'}_{S,\chi}\br{s-\rho,\tau\otimes\sigma}}{Z_{S,\chi}\br{s-\rho,\tau\otimes\sigma}}=O\br{t^{n-1}}
\end{equation}
\newline
for $s=\sigma_{1}+\II t$, $\rho+u\leq\sigma_{1}<\frac{1}{2}t+\rho$.

Finally, by (\ref{3.1}), (\ref{4.12}) and (\ref{4.13})

\begin{equation}\label{4.14}
\frac{Z^{'}_{R,\chi}\br{s,\sigma}}{Z_{R,\chi}\br{s,\sigma}}=\sum\limits_{p=0}^{n-1}\br{-1}^{p}\sum\limits_{\br{\tau,\lambda}\in I_{p}}\frac{Z^{'}_{S,\chi}\br{s+\rho-\lambda,\tau\otimes\sigma}}{Z_{S,\chi}\br{s+\rho-\lambda,\tau\otimes\sigma}}=
\end{equation}
\[\sum\limits_{p=0}^{n-1}\br{-1}^{p}\sum\limits_{\substack{\br{\tau,\lambda}\in
I_{p}\\\lambda<2\rho}}\frac{Z^{'}_{S,\chi}\br{s+\rho-\lambda,\tau\otimes\sigma}}{Z_{S,\chi}\br{s+\rho-\lambda,\tau\otimes\sigma}}+\]
\[\sum\limits_{p=0}^{n-1}\br{-1}^{p}\sum\limits_{\br{\tau,2\rho}\in I_{p}}\frac{Z^{'}_{S,\chi}\br{s-\rho,\tau\otimes\sigma}}{Z_{S,\chi}\br{s-\rho,\tau\otimes\sigma}}=\]
\[\sum\limits_{p=0}^{n-1}\br{-1}^{p}\sum\limits_{\substack{\br{\tau,\lambda}\in
I_{p}\\\lambda<2\rho}}O\br{t^{n-1}}+\]
\[\sum\limits_{p=0}^{n-1}\br{-1}^{p}\sum\limits_{\br{\tau,2\rho}\in I_{p}}O\br{t^{n-1}}+
\sum\limits_{p=0}^{n-1}\br{-1}^{p}\sum\limits_{\br{\tau,2\rho}\in I_{p}}\sum\limits_{\abs{t-\gamma_{S,p,\tau}}\leq 1}\frac{1}{s-\rho_{S,p,\tau}}=\]
\[O\br{t^{n-1}}+\sum_{\abs{t-t_{R}}\leq 1}\frac{1}{s-s_{R}}\]
\newline
for $s=\sigma_{1}+\II t$, $\rho\leq\sigma_{1}<\frac{1}{2}t+\rho$. This proves (\textit{a}).

(\textit{b}) Let $u>0$. Obviously,

\[\abs{\sum\limits_{\abs{t-t_{R}}\leq 1}\frac{1}{s-s_{R}}}\leq\sum\limits_{\abs{t-t_{R}}\leq 1}\frac{1}{\abs{s-s_{R}}}<\frac{1}{u}\sum\limits_{\abs{t-t_{R}}\leq 1}1\]
\newline
for $s=\sigma_{1}+\II t$, $\rho+u\leq\sigma_{1}<\frac{1}{2}t+\rho$. Hence, it follows from (\ref{3.5}) and (\ref{4.14}) that
\[\frac{Z^{'}_{R,\chi}\br{s,\sigma}}{Z_{R,\chi}\br{s,\sigma}}=O\br{t^{n-1}}\]
\newline
for $s=\sigma_{1}+\II t$, $\rho+u\leq\sigma_{1}<\frac{1}{2}t+\rho$. This completes the proof.
\end{proof}

\begin{remm}{4.2.}\label{rem}
Approximate formulas for the logarithmic derivative of the zeta functions were quite often exploited by many authors (see, e.g., \cite{Hejhal1}-\cite{Titchmarsh}), not always for the same underlaying space. Usually, they were applied to obtain error terms in the prime number resp. prime geodesic theorem, where the search for the optimal error bound is widely open (see, e.g., \cite{AG0,Park}).
\end{remm}

\end{document}